\documentclass{ws-ijbc}
\usepackage{ws-rotating}     
\usepackage{graphicx}
\usepackage{epstopdf}
\usepackage{color}
\usepackage{multirow}
\usepackage{amsmath}
\begin{document}
\interfootnotelinepenalty=10000
\newcommand{\LvV}[1]{{\color{red} #1}}
\newcommand{\MH}[1]{{\color{green} #1}}

\catchline{}{}{}{}{} 

\markboth{Hoti \& van Veen {\sl et al.}}{Saddle-node--transcritical interactions}

\title{Saddle-node--transcritical interactions in a stressed predator-prey-nutrient system}

\author{Lennaert van Veen}

\address{Faculty of Science, University of Ontario Institute of Technology\\2000 Simcoe Street North, L1H7K4 Oshawa, ON, Canada}

\author{Marvin Hoti}

\address{Department of Mathematics, Ryerson University\\350 Victoria Street, M5B 2K3 Toronto, ON, Canada\\marvin.hoti@ryerson.ca}

\maketitle

\begin{history}
\received{(to be inserted by publisher)}
\end{history}

\begin{abstract}
We examine the interaction of transcritical and saddle-node bifurcations in a predator-prey-nutrient system that is stressed by the presence of a toxicant
affecting the prey. This model, formulated by Kooi et al. ({\sl Ecol. Model.} {\bf 212}(2008), 304--318), has a two-dimensional invariant sub system with zero predator
density. In the sub system, a pair of prey-nutrient equilibria is created in a saddle-node bifurcation, while predator invasion in modelled by a transcritical bifurcation
of one of this pair. Interactions of these bifurcations at codimension-two points give rise to bistable, periodic and heteroclinic predator-prey-nutrient dynamics. We explain why
the the codimension-two points are numerically detected as cusp and Bogdanov-Takens points when using standard test functions and propose a new test function for systems 
with codimension one trancritical curves.
\end{abstract}

\keywords{Predator-prey-nutrient-toxicant dynamics, saddle-node--transcritical bifurcation, test functions, MatCont, AUTO.}

\section{Introduction}\label{intro}
\noindent Dynamical systems with invariant axes and (hyper) planes arise in various classes of models, such as in disease models (e.g. the malaria models reviewed by \citet{Teboh2013}), mass action kinetics (going all the way back to \citet{Lotka1920}) and plasma physics (see, e.g., \citet{Bian2003}). The invariance reflects an elementary property of the model. For instance, if no infected individuals are present, the disease cannot spread; a certain chemical complex can be recycled but not created in a reaction chain and turbulent fluctuations grow through self-interaction. A consequence of this special structure is that transcritical bifurcations can occur when varying a single parameter, just like saddle-node bifurcations in generic dynamical systems. 

In the context of population dynamics, the transcritical bifurcation is sometimes called {\em invasion}, as it marks the onset of a ``positive'' equilibrium, i.e. an equilibrium with positive values of the unknowns, which models co-existence of species\footnote{In this paper, a positive (negative, non-negative) solution (e.g. equilibrium or periodic) is a solution for which all variables remain positive (negative, non-negative) for all time.}. In the other contexts mentioned above, this would correspond to co-existence of infected and healthy individuals, multiple chemical complexes or background shear flow and turbulent fluctuations. 

\begin{figure}[t]
\begin{center}
\includegraphics[width=250pt]{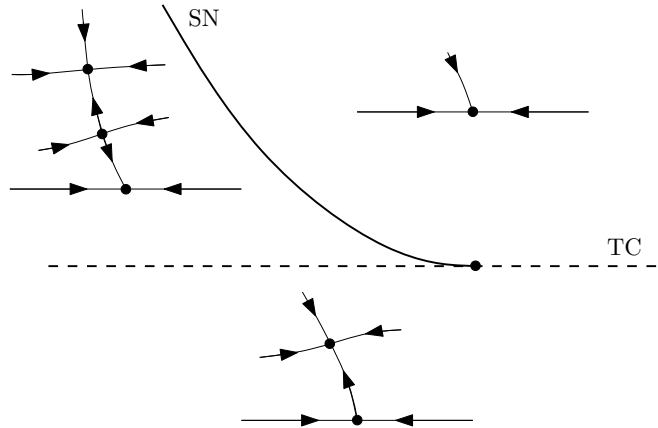}
\end{center}
\caption{Unfolding of the single zero SNTC bifurcation including only the non-negative equilibria. In the bifurcation diagram, the dashed line represents the transcritical bifurcation and the solid line represents the saddle-node bifurcation. In the phase portraits, the solid horizontal line represents an invariant axis. Due to the similarity to the unfolding of the generalised Hopf bifurcation, this singularity has been labeled ``generalised transcritical bifurcation''. In the complete diagram, found in \citet{Saputra2010}, the saddle-node bifurcation does not terminate. Note, that this singularity can occur in a system with a single degree of freedom. We have drawn two-dimensional phase portraits for easy comparison to the double zero case in Figs. \ref{unfolding1}-\ref{unfolding3}.}
\label{single_zero}
\end{figure}

Following the general philosophy of bifurcation analysis, it is natural to ask the question what interactions the transcritical bifurcation can have with other singularities. Such interaction points can act as ``organising centres'', tying together various codimension one bifurcations and organising the qualitative dynamics for ranges of parameter values. The simplest possibility is to have a zero eigenvalue along with a vanishing normal form coefficient. At such a codimension two point, curves of transcritical bifurcations and saddle-node bifurcations are tangent. This case can be realised in a model with a single degree of freedom and was investigated in detail by \citet{Saputra2010} and \citet{Saputra2015}. When one only considers non-negative equilibria, the unfolding looks similar to that of the generalised Hopf bifurcation, which explains why this singularity is sometimes called the generalised transcritical bifurcation. This is demonstrated in Fig. \ref{single_zero}. Here, we will refer to it as the single zero SNTC interaction.

\begin{figure}[t]
\begin{center}
\includegraphics[width=300pt]{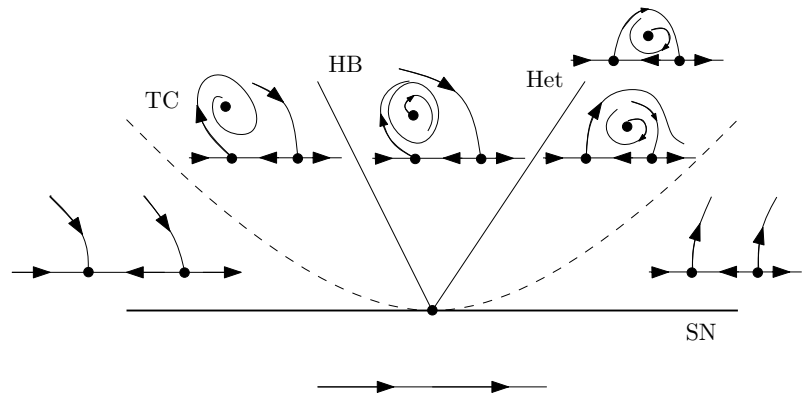}
\end{center}
\caption{Unfolding of the SNTC interaction with a double zero eigenvalue including only the non-negative equilibria. Shown is the saddle case with a positive periodic orbit.
The bifurcations in the diagrams have been labelled SN for saddle-node, TC for transcritical, HB for Hopf and Het for heteroclinic. In the phase portraits, the solid horizontal line represents an invariant axis. }
\label{unfolding1}
\end{figure}

\begin{figure}[t]
\begin{center}
\includegraphics[width=300pt]{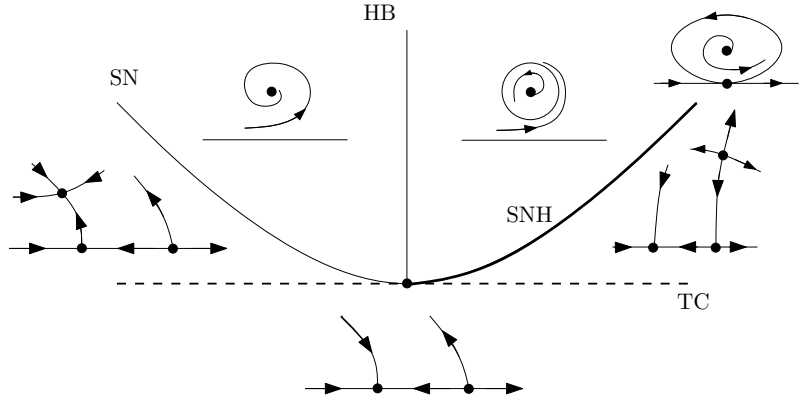}
\end{center}
\caption{Unfolding of the SNTC interaction with a double zero eigenvalue including only the non-negative equilibria. Shown is the elliptic case with a positive periodic orbit.
The bifurcations in the diagrams have been labelled SN for saddle-node, TC for transcritical, HB for Hopf and SNH for a homoclinic to a saddle-node. In the phase portraits, the solid horizontal line represents an invariant axis. }
\label{unfolding2}
\end{figure}

\begin{figure}[t]
\begin{center}
\includegraphics[width=300pt]{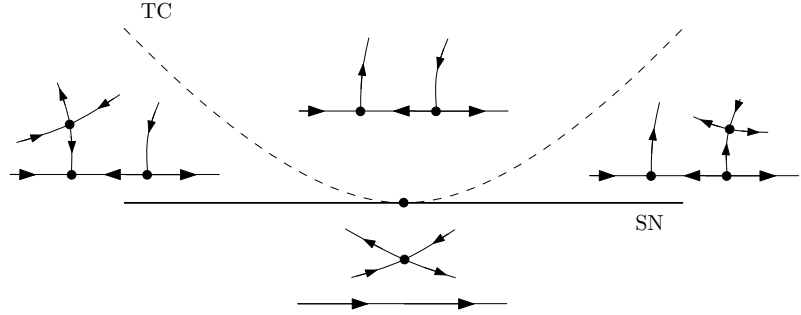}
\end{center}
\caption{Unfolding of the SNTC bifurcation with a double zero eigenvalue including only the non-negative equilibria. Shown is the saddle or elliptic case with a negative periodic orbit.
The bifurcations in the diagrams have been labelled SN for saddle-node and TC for transcritical. In the phase portraits, the solid horizontal line represents an invariant axis. }
\label{unfolding3}
\end{figure}

A more complicated interaction requires at least two degrees of freedom and involves a zero eigenvalue with algebraic multiplicity two and geometric multiplicity one. The corresponding eigenvector lies in the invariant plane, while the generalised eigenvector is transversal to it. Normal forms and unfoldings for the double zero SNTC interaction were presented by \citet{Saputra2010} and \citet{Saputra2015} and involve, in addition to the saddle-node and transcritical bifurcations, a Hopf bifurcation and a curve along which a periodic orbit disappears. The latter can be either a heteroclinic bifurcation or a homoclinic to a saddle-node. When we omit the negative equilibria, we obtain one of the three unfoldings shown in Figs. \ref{unfolding1}-\ref{unfolding3}.

The third possibility, which can occur in systems of dimension three or more, involves a single zero and a pair of purely imaginary eigenvalues. This transcritical-Hopf case was investigated in some detail in the literature, starting from \citet{Langford1979}, who presented several one-parameter unfoldings near the singularity. Later, \citet{Jiang2010} and \citet{Saputra2015} presented four distinct two-parameter unfoldings. In the vicinity of this point, one can find periodic orbits inside and outside the invariant plane, as well as quasi-periodic motion.

Out of these three cases, the first and the last have been reported on in model studies. The first requires only a single boundary equilibrium and two positive (or negative) equilibria.
It was found to play a role, for instance, in a model of smoking as an epidemic by \citet{Voorn2013} and one of dispersal patterns for co-existing species by \citet{Mohd2018}. The the latter case, the singularity is presented as a ``triple point'' which appears on the mutual boundary of regions of extinction, co-existence and bi-stability between these modes. It is also present in the model for plasma physics mentioned before and, in fact, the single-degree-of-freedom model formulated by \citet{Bian2003} is identical to the normal form proposed by \citet{Saputra2010}.

To mention but a few examples of the analysis of the transcritical-Hopf bifurcation: \citet{Langford1979} found it in a model for fluid motion due to Hopf, \citet{Doedel1984} in a predator-prey-nutrient model and \citet{Gimmelli2015} in an ecoepidemic model. 
It was also identified in the stressed predator-prey-nutrient model under consideration here. \citet{Kooi2008} showed that, around this point, stable equilibrium and periodic solutions exist with zero or positive predator density. From the layout of the Hopf and transcritical bifurcations, we can tell that the unfolding corresponds to case III of \citet{Saputra2015}.

To the best of our knowledge, the SNTC case with two zero eigenvalues has never been identified in an actual application. One may speculate the this is in part due to the fact that it is neither as easy to analyse as the case with a single zero eigenvalue nor as much part of the canon of singularity theory as the case with complex conjugate eigenvalues, having been analysed relatively recently. We show that a change in a single model parameter, related to the growth rate of predators, introduces this point as an organising centre in the model by \citet{Kooi2008}. In its vicinity, there exists a limit cycle of arbitrarily large period that models switching between two nutrient-prey equilibria, mediated by the predator population. 

In addition to presenting the first analysis of the double zero SNTC interaction in an application, a secondary goal of this paper is to explain why the widely used bifurcation analysis software packages MatCont \cite{Dhooge2004} and AUTO \cite{AUTO} erroneously classify the single-zero and double-zero SNTC cases as cusp and Bogdanov-Takens (BT) points, respectively. We suggest a test function that takes the special structure of the model into account and uniquley identifies the SNTC interactions.

\section{The stressed model}\label{model}
We consider the stressed model introduced by \citet{Kooi2008}. Here, "stressed" refers to the presence of a toxicant in the water flowing into a basin that contains a nutrient as well as a prey and a predator population. The basin is considered to be well-mixed and hold a constant volume of water. The densities of the nutrient, prey, predator and toxicant then satisfy
\begin{align}
\frac{dN}{dt}&=(N_r-N)D-I_{NR}\frac{N}{\kappa_{NR}+N}R, \label{dNdt}\\
\frac{dR}{dt}&=\bigg(\mu_{NR}\frac{N}{\kappa_{NR}+N}-(D+m_R(c_R))\bigg)R-I_{RP}\frac{R}{\kappa_{RP}+R}P, \label{dRdt}\\
\frac{dP}{dt}&=\bigg(\mu_{RP}\frac{R}{\kappa_{RP}+R}-(D+m_{P0})\bigg)P,\label{dPdt}\\
\frac{dc_{T}}{dt}&=(c_{r}-c_{T})D, \label{dcTdt}
\end{align}
respectively. The nutrient density tends to relax to the inflow value, $N_r$, which is the first control parameter, at the flow rate, $D$, which is the second control parameter. It is also consumed by the prey which, in turn, is consumed by the predator. The consumption of the nutrient by the prey and the prey by the predator is modelled by a Holling type II functional response. Both the prey and the predator are drained from the basin a the flow rate $D$, and die at a rate that is fixed for the predator but depends on the toxicant concentration for the prey. The total toxicant can be split up into the portion ingested by the prey, $c_R$, and the portion in the ambient water, $c_W$, according to
\begin{equation}
 c_{T}=c_{W}+c_R R. \label{cTcWcR}
\end{equation}
On time scales relevant to this model, the proportion of toxicant carried by the prey is fixed by 
\begin{equation} 
c_R=BCF_{WR} \,c_W.   \label{cRcW}
\end{equation}
Finally, the maintenance rate constant of the prey is given by
\begin{equation}
m_{R}(c_{R})=m_{R0}\bigg(1+\frac{\max\{c_{R}-c_{RM0},0\}}{c_{RM}} \bigg). \label{mR}
\end{equation}
Together, Eqs. (\ref{dNdt}--\ref{mR}) form a closed system. However, since we are interested in equilibria, limit cycles and connecting orbits rather than transient motion, we will assume that $c_T$ has assumed its equilibrium value $c_r$ and eliminate the toxicant concentrations to find
\begin{align}
\frac{dN}{dt}&=(N_r-N)D-I_{NR}\frac{N}{\kappa_{NR}+N}R, \label{ODE1} \\
\frac{dR}{dt}&=\bigg(\mu_{NR}\frac{N}{\kappa_{NR}+N}-(D+m_R(R))\bigg)R-I_{RP}\frac{R}{\kappa_{RP}+R}P, \label{ODE2}\\
\frac{dP}{dt}&=\bigg(\mu_{RP}\frac{R}{\kappa_{RP}+R}-(D+m_{P0})\bigg)P, \label{ODE3}
\end{align}
where
\begin{equation}
\label{mRR}
m_R(R)=m_{R0}+\frac{m_{R0}}{c_{RM}}\max\left\{\frac{BCF_{CW} c_r}{1+BCF_{CW} R}-c_{RM0},0\right\}.
\end{equation}
Table \ref{pars} lists the parameters of the model with their interpretation, value and dimension.
\begin{table}\label{pars}
\tbl{Parameter set for a stressed bacterium-ciliate model. All values are taken from \citet{Kooi2008} except for $\mu_{RP}$, which is changed to locate the double-zero SNTC interaction. The units are given in terms of time, $t$, the mass of toxicant, $m$, the volume of the basin, $\nu$, and the (bio)volume of the nutrient, predator and prey, $V$.}
{\begin{tabular}{|llll|}\hline
$\mu_{NR}$ & Max growth rate & $t^{-1}$ & $0.5$ $h^{-1}$\\
$I_{NR}$ & Max ingestion rate & $t^{-1}$ & $1.25$ $h^{-1}$\\
$\kappa_{NR}$ & Saturation constant & $V/\nu$ & $8.0$ $mg/dm^{3}$\\
$m_{R0}$ & Maintenance rate coefficient & $t^{-1}$ & $0.025$ $h^{-1}$\\
$c_{RM0}$ & No effect concentration & $m/V$ & $0.1$ $\mu g/mg$\\
$c_{RM}$ & Tolerance concentration  & $m/V$ & $0.5$ $\mu g/mg$\\
$BCF_{WR}$ & Bioconcentration factor & $\nu /V$ & $1.0$ $dm^{3}/mg$ \\
$\mu_{RP}$ & Max growth rate & $t^{-1}$ & $0.2$ or $0.7$ $h^{-1}$\\
$I_{RP}$ & Max ingestion rate & $t^{-1}$ & $0.333$ $h^{-1}$\\
$\kappa_{RP}$ & Saturation constant & $V/\nu$ & $9.0$ $mg/dm$\\
$m_{P0}$ & Maintenance rate coefficient & $t^{-1}$ & $0.01$ $h^{-1}$\\
$c_r$ & Toxicant concentration at inflow & $m/\nu$ & $1$ or $9$ $\mu g/dm^3$\\
\hline 
\end{tabular}}
\end{table}

Two invariant subsystems exist in this model. In the invariant subsystem $R=P=0$ all solutions tend towards the ``wash-out'' equilibrium $(N_r, 0, 0)$. 
The invariant predator-free subsystem $P=0$ can support two equilibria in addition to the wash-out solution, as well as limit cycles and connecting
orbits. The full system can support equilibria, cycles and connecting orbits with co-existing predator and prey populations.


\section{The single zero SNTC interaction}

The first organising centre identified by \citet{Kooi2008} is a single zero SNTC interaction in the nutrient-prey system at $c_r=1$ and $\mu_{RP}=0.2$. In their Fig. 2, it is labeled ``N'' and the unfolding is that shown in our Fig. \ref{single_zero}. The three regions around the codimension two point are characterised by wash-out, a stable positive prey population and co-existence of these stable solutions.

\subsection{Test functions for the single zero SNTC interaction}

When using MatCont \cite{Dhooge2004} to construct the bifurcation diagram, the codimension two point is classified as cusp point on the saddle-node curve. Here, we verify that the test function used by MatCont has an isolated zero at the SNTC point and we suggest a test function to distinguish the two cases.

Suppose we are tracing a saddle-node bifurcation in two parameters, $\lambda_1$ and $\lambda_2$ in a system with $n$ variables. Along the curve we have
\begin{equation}
f(x,\lambda_1,\lambda_2)=0,\ D_x f(x, \lambda_1,\lambda_2)q=0,
\label{defSN}
\end{equation}
where $q$ is the right zero eigenvector. In practice, the left and right zero eigenvalue are computed by solving the following nonsingular bordered systems \cite{Govaerts2000}
\begin{align}\label{bordered1}
\begin{bmatrix}
p^t &\gamma
\end{bmatrix}\begin{bmatrix}
Df(x,\lambda_1,\lambda_2) & \bar{p}\\
\bar{q}^t & 0
\end{bmatrix}
&=
\begin{bmatrix}
0 &1
\end{bmatrix},
&\quad
\begin{bmatrix}
Df(x,\lambda_1,\lambda_2) & \bar{p}\\
\bar{q}^t & 0
\end{bmatrix}
\begin{bmatrix}
q\\\gamma
\end{bmatrix}&=
\begin{bmatrix}
0\\1
\end{bmatrix}.
\end{align} 
The eigenvectors are thus normed as $\bar{p}^t p=\bar{q}^t q=1$, where $\bar{p}$ and $\bar{q}$ are the zero eigenvectors computed at the previous point on the saddle-node curve.
In a similar fashion, a vector tangent to the solution curve is found. We will denote it by $v$ and assume it is normalised to unit length. Finally, the Hessian of the dynamical system is computed. Following \citet{Kuznetsov1999}, we will consider it as a bilinear form $\mathbb{R}^n\times\mathbb{R}^n\rightarrow \mathbb{R}^n$ denoted by $D^2_{x,x}f$.

Generically, the dynamics in the centre manifold is conjugate to  
the normal form
\begin{align}
\dot{y}&= \alpha a+\beta y^2\\
\alpha&=\langle p, \left. D_s f(x,\lambda_1(s),\lambda_2(s))\right|_{s=0}\rangle\\
\beta&=\langle p, (D_{x,x}f)(q,q)\rangle
\end{align}
where the brackets denote the dot product and $D^2_{x,x}f$ is evaluated at $s=0$, where the saddle-node bifurcation takes place. 

If $\beta=0$ at some point along the curve, MatCont classifies it as a cusp. However, if $\alpha=0$ at the same point, it is actually a SNTC point. Thus, we propose to use $\alpha$ as additional test function. A suitable choice for the unfolding curve is one perpendicular to the saddle-node curve, e.g. $(s v_2,-s v_1)$. The test function then evaluates to
\begin{equation}\label{TF}
\alpha=\langle p, D_{\lambda_1}f v_2-D_{\lambda_2}f v_1\rangle
\end{equation}

The unfolding of the cusp bifurcation is given by
\begin{equation}
\dot{x}=a + b x +x^3.
\end{equation}
Simply taking $p=q=1$, we find along the saddle-node curve, defined by $a^2/4+b^3/27=0$, that $\beta=2^{2/3} a^{1/3}/9$ while $\alpha$ stays bounded away from zero. In fact, $\alpha=1$ at the bifurcation point.

The normal form for the single zero SNTC interaction proposed by \citet{Saputra2010} and \citet{Saputra2015} is
\begin{equation}
\dot{x}=a x+ b x^2 +x^3.
\end{equation}
The equilibrium at $x=0$ undergoes a transcritical bifurcation along $a=0$ and a saddle-node bifurcation along the line $a=b^2/4$.
Along the latter curve, we find $\alpha=-b/4+\mbox{O}(b^3)$ and $\beta=-b$ so both test functions are zero at the SNTC point.

Should one want to detect the SNTC interaction along a transcritical bifurcation curve, one should monitor test function $\beta$. By default, no test functions are computed along transcritical curves in MatCont.

\section{The double zero SNTC interaction}

\begin{figure}[t]
\label{bifcurve}
\begin{center}
\includegraphics[width=300pt]{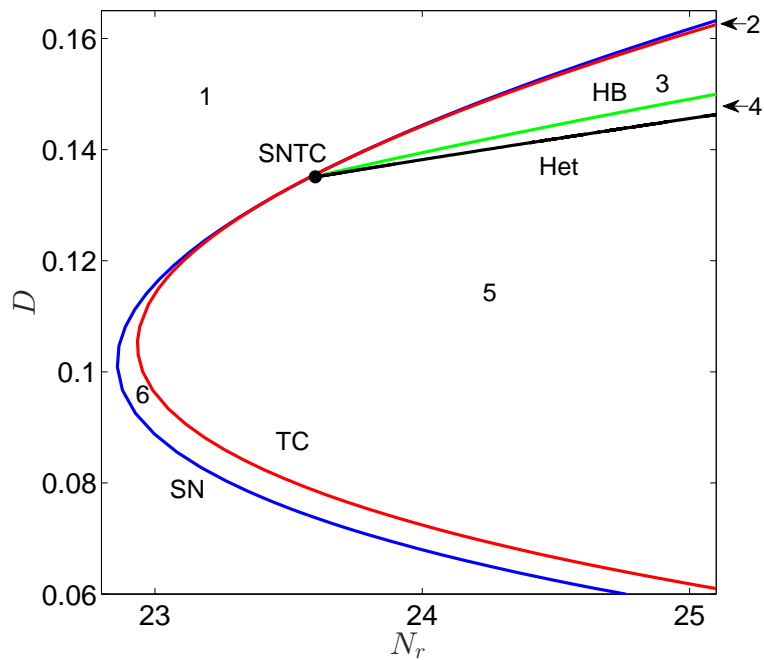}
\end{center}
\caption{Partial bifurcation diagram of the stressed nutrient predator prey system with $N$ and $D$ as the free parameters. The diagram was computed using MatCont \cite{Dhooge2004} at $c_r=9$ and $\mu_{RP}=0.7$. The bifurcations have been labeled as in Fig. \ref{unfolding1}.}
\label{2zero}
\end{figure}
The SNTC interaction with two zero eigenvalues occurs in the stressed model at parameter values $c_r=9$ and $\mu_{RP}=0.7$. Figure \ref{2zero} shows a partial bifurcation diagram in the nutrient density of the inflow, $N_r$, and the flow rate, $D$. It follows the pattern of unfolding (\ref{unfolding1}). The wash-out equilibrium is stable in this entire diagram. In region 1, this is the only equilibrium. In region 2, two equilibria with nonzero prey density exist. The one with the highest prey density is stable, while the other is of saddle type. The stable manifold of the latter separates the basins of attraction of the stable wash-out and nonzero prey equilibria. In region 3, the predator has invaded the system, rendering both equilibria on the boundary $P=0$, $R>0$ unstable. 

\begin{figure}[t]
\label{timeseries}
\begin{center}
\includegraphics[width=270pt]{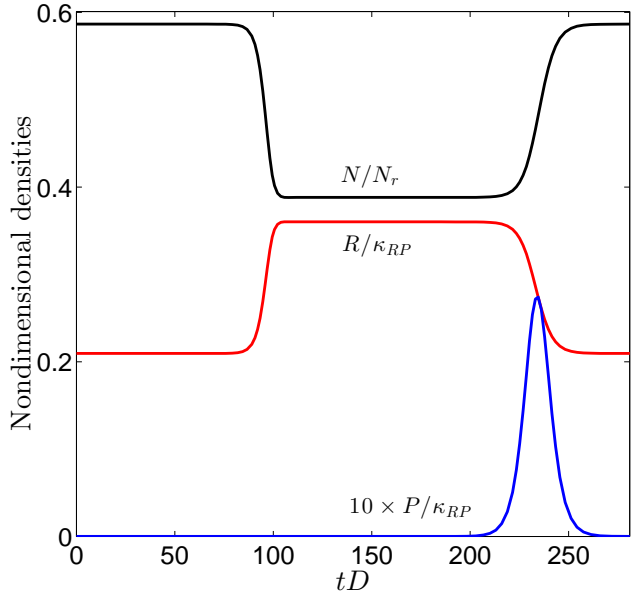}
\end{center}
\caption{Time series in region 4 of diagram \ref{2zero}, near the heteroclinic bifurcation curve. Shown are the nutrient, normalised by $N_r$, in black and the prey and predator densities, normalised by $\kappa_{RP}$, in red and blue, respectively. The parameter values are $N_r=24.6$ and $D=0.142$.}
\label{ts}
\end{figure}

In region 4, the co-existence equilibrium has turned unstable in a Hopf bifurcation and the system periodically oscillates. Close to the heteroclinic bifurcation, the period of the oscillation grows without bound, and the observed behaviour looks like switching between the equiilbria with high and low prey density in the boundary. An example is shown in Fig. \ref{ts}. In the first phase of the switching, the prey density grows while the predator density stays close to zero. In the next phase, the predator density slowly grows while the prey density remains nearly constant and high. The predators then thrive and rapidly consume the prey until the density of the latter drops below a critical threshold, after which the predator population declines. 

In region 5, the cycle has disappeared in a heteroclinic bifurcation. After this bifurcation, the stable manifold of the zero predator equilibirum with the higher prey density no longer bounds the basin of attraction of the wash-out solution, which seems to attract all initial conditions. This remains true in region 6, where the coexistence equilibrium is no longer positive. The main difference between regions 5 and 6, on one hand, and 1, on the other, is that in regions 5 and 6 long transients can be observed due to the existence of unstable equilibria with nonzero predator in prey densities. 

\subsection{Test functions for the double zero case}
The test function used by MatCont to locate BT points along saddle-node curves is extracted from the following nonsingular bordered systems
\begin{align}\label{bordered2}
\begin{bmatrix}
u^t &\beta
\end{bmatrix}\begin{bmatrix}
Df(x,\lambda_1,\lambda_2) & \bar{p}\\
\bar{q}^t & 0
\end{bmatrix}
&=
\begin{bmatrix}
p^t &0
\end{bmatrix},
&\quad
\begin{bmatrix}
Df(x,\lambda_1,\lambda_2) & \bar{p}\\
\bar{q}^t & 0
\end{bmatrix}
\begin{bmatrix}
v\\\beta
\end{bmatrix}&=
\begin{bmatrix}
q\\0
\end{bmatrix},
\end{align} 
that are solved after systems (\ref{bordered1}). If, and only if, $\beta=0$, $u$ and $v$ are left and right generalized eigenvectors for eigenvalue zero. Multiplying the second system from the left by $(p^t,0)$, we find that $\beta=p^t q$ so that the condition $\beta=0$ can also be interpreted as the solvability condition for the existence of such generalised eigenvectors. A BT point is detected when $\beta=0$.

Two models have been proposed for the double zero SNTC interaction. In the absence of a formal derivation, \citet{Saputra2010} called the following system a ``minimal model'':
\begin{align}
\dot{x}&= y \nonumber\\
\dot{y}&= a x + k_1 b y + b x^2 + k_2 x y+x^2 y+\epsilon x^3+k_3 x^4, \label{SNTC1}
\end{align}
where the constants $k_1,\, k_2,\, k_3\neq 0$ satisfy
\begin{alignat}{2}
2\epsilon k_1^2-k_1 k_2&=1,  \nonumber\\
3 k_1 k_3 &=1. \label{conditions}
\end{alignat}
In this model, the SNTC interaction takes place at $a=b=0$ and the saddle-node bifurcation takes place along the line $a=\epsilon (2\sqrt{(1-3 k_3 b)^3} -2+9 k_3 b)/(27 k_3^2)$. A straightforward computation shows that, along this curve, $\beta=\beta_1\epsilon k_3 b+\mbox{O}(b^2)$, where $\beta_1$ is a constant that depends on the choice of eigenvectors, which explains why the SNTC interaction is classified as a BT point. A more formally derived normal form was presented by \citet{Saputra2015}:
\begin{align}
\dot{x}&=a +y+ \mu_1 x^2 \nonumber\\
\dot{y}&=b y +\mu_2 x y. \label{SNTC2}
\end{align}
Again, the SNTC interaction takes place at $a=b=0$, and the saddle-node bifurcation takes place along the curve $a=0$. It is easy to see that, along this curve, $\beta=\beta_1 b+\mbox{O}(b^2)$. Again, the BT test function has a zero at the SNTC point.

At a zero of the test function $\beta$, MatCont automatically computes the coefficient $a_2$ of the BT normal form
\begin{align}
\dot{x}&=y  \nonumber\\
\dot{y}&=\lambda_1+\lambda_2 x+a_2 x^2+b_2 x y. \label{BT}
\end{align}
As is apparent from model (\ref{SNTC1}) of the SNTC interaction, $a_2$ evaluates to zero at the codimension two point. Thus, $a_2$ could be used as an additional test function. However, this coefficient is also zero in the case of a degenerate BT bifurcation of codimension three, which does not involve a transcritical bifurcation. It is better to use test function (\ref{TF}), that relies on the absence of a constant term in the unfolding. 

A direct computation shows that, for SNTC models (\ref{SNTC1}) and (\ref{SNTC2}), $\alpha=\alpha_1 \epsilon b+\mbox{O}(b^2)$ and $\alpha=\alpha_1 b+\mbox{O}(b^2)$, repspectively, where $\alpha_1$ is a numerical constant that depends on the choice of eigenvectors. In contrast, for the BT normal form, we find $\alpha=\|p\|$ at the SNTC point. Thus, test function (\ref{TF}) neatly distiguishes the SNTC interaction from the BT bifurcation.

When tracing out the transcritical bifurcation curve, one can simple use $\beta$ as a test function to localize the double zero SNTC interaction.

\section{Conclusion}

We have analysed the single zero and double zero SNTC interaction in the stressed nutrient-predator-prey model of \citet{Kooi2008}, and demonstrated that both play the role of ``organsing centre'' for the dynamics. Since the transcritical bifurcation is of codimension one in this model, due to its special structure, the numerical construction of the bifurcation diagram is not entirely straightforward. The 
single zero SNTC interaction is classified by MatCont \cite{Dhooge2004} as a cusp bifurcation, while the double zero interaction is classified as a BT point. We verified that the applicable test functions have isolated zeros for the normal forms of the SNTC interactions proposed by \citet{Saputra2015} and \citet{Saputra2010}. While we have refered to MatCont in the discussion, AUTO \cite{AUTO} uses the same test functions, albeit computed in a slightly different way, aproximating $\beta$ by a finite-difference formula. Our observations about the test functions are thus valid for AUTO, too.

We propose a new test function that can uniquely identify the SNTC interactions along saddle-node curves. This test function only involves the parametric derivatives of the dynamical system and the null space of its Jacobian. We hope that this test function will prove a useful tool for performing numerical bifurcation analysis on models with invariant coordinate planes, often found in the modelling of population dynamics and chemical reactions.


\nonumsection{Acknowledgements} \noindent LvV was supported by an individual Discovery Grant of NSERC.

\end{document}